\documentclass[leqno]{amsproc}
\usepackage{amsfonts}
\usepackage{amssymb, latexsym, amsmath, pb-diagram}
\usepackage{pstricks-add}
\usepackage{multicol}
\usepackage{bm}
\usepackage[mathscr]{eucal}
\usepackage{xy}
\usepackage{epic,eepic}
\xyoption{all}

\def\w{\widetilde}

\def\o{\overline}

\newcommand{\zz}{\mathbb{Z}}

\newcommand{\RR}{\mathbb{R}}

\begin{document}
\title[Simplicial Homeology and Homeotopy]
{\bf Simplicial Homeology and Homeotopy}
\author[Q. Zheng, F. Fan]{Qibing Zheng and Feifei Fan}
\thanks{The author is supported by NSFC grant No. 11071125}
\date{June 19, 2011}
\address{School of Mathematical Science and LPMC, Nankai University,
Tianjin, 300071, P.R.China}
\email{zhengqb@nankai.edu.cn}
\maketitle
\vspace{-6mm}
\begin{abstract} In this paper, we define homeology group, reduced homeology group, cohomeology group and reduced cohomeology group on finite simpicial complexes and prove that these groups are homeomorphism invariants of polyhedra. We also define homeotopy type of polyhedra which is finer than homotopy type but coarser than homeomorphism class.
\end{abstract}\vspace{3mm}

Simplicial theory are getting more and more important as the toric topology develops. But so far, there are only algebraic homotopy invariants of polyhedra. The torsion group of face ring is a simplicial invariant that is not a homeomorphism invariant. In this paper, we define homeology  group $H_{s,t}^r$, reduced homeology  group $\w H_{s,t}^r$, cohomeology group $H^{s,t}_r$ and reduced cohomeology group $\w H^{s,t}_r$ for all simplicial complexes and prove that these groups are homeotopy invariants of polyhedra that are not homotopy invariants but are homeomorphism invariants. In fact, (co)homology is the limit (the convergence group) of (co)homeology (spectral sequence). An equivalent homeotopic relation is defined between solid maps and homeotopy type of polyhedra is  defined on homeotopy. Homeotopy types of polyhedra are finer than homotopy types but coarser than homeomorphism classes. For example, (co)homeology group distinguishes disks of different dimension and many other contractible polyhedra.

\section{Homeology and cohomeology spectral sequence}

The following definition is a brief review of finite simplicial complex theory.\vspace{3mm}

{\bf Definition 1.1} A (finite, abstract) simplicial complex $K$ with vertex set $S$
is a set of subsets of the finite set $S$ satisfying the following two conditions.

1) For all $s\in S$, $\{s\}\in K$.

2) If $\sigma\in K$ and $\tau\subset\sigma$, then $\tau\in K$. Specifically, the empty set $\phi\in K$.

A simplicial subcomplex $L$ of $K$ is a subset $L\subset K$ such that $L$ is also a simplicial complex.

For simplicial complexes $K$ and $L$ with respectively vertex set $S$ and $T$, a simplicial map $f\colon K\to L$ is a map $f\colon S\to T$ such that for all $\sigma\in K$, $f(\sigma)\in L$. $K$ is simplicial isomorphic to $L$ if there are simplicial maps $f\colon K\to L$ and $g\colon L\to K$ such that $fg=1_L$ and $gf=1_K$ ($1_K$ and $1_L$ denote the identity map).

An element $\sigma$ of a simplicial complex $K$ is called an $n$-simplex of $K$ if $\sigma$ has $n{+}1$ vertices of $S$ and $|\sigma|=n$ is the dimension of $\sigma$. Specifically, the empty simplex $\phi$ has dimension $-1$. The dimension of $K$ is the biggest dimension of its simplexes.  A face of a simplex $\sigma$ is a subset $\tau\subset\sigma$. A proper face is a proper subset. For a simplex $\sigma$ of $K$, the link of $\sigma$ is the simplicial subcomplex link$_{K}\sigma=\{\tau\in K\,|\,\sigma{\cup}\tau\in K,\,\sigma{\cap}\tau=\phi\}$ and $2^{\sigma}$ is the simplicial subcomplex of $K$ consisting of all faces of $\sigma$ and $\partial\sigma$ is the  simplicial subcomplex of $K$ consisting of all proper faces of $\sigma$.

The geometrical realization of a simplicial complex $K$ is the topological space $|K|$ defined as follows. Suppose the vertex set $S=\{v_1,\cdots,v_m\}$. Let $e_i\in\RR^m$ be such that the $i$-th coordinate of $e_i$ is $1$ and all other coordinates of $e_i$ are $0$. $|K|$ is the union of all the convex hull of $e_{i_0},\cdots,e_{i_s}$ such that $\{v_{i_0},\cdots,v_{i_s}\}$ is a simplex of $K$. The geometrical realization of the simplicial map $f\colon K\to L$ is the map from $|K|$ to $|L|$ that linearly extends $f$ which is still denoted by $f$. The convex hull of $e_{i_0},\cdots,e_{i_s}$ is a geometrical $s$-simplex of $|K|$. A polyhedron $X$ is a topological space that is homeomorphic to the geometrical realization $|K|$ of a simplicial complex $K$ and $K$ is a triangulation of $X$.

A stellar subdivision of $K$ on the simplex $\sigma$ is the simplicial complex $S_{\sigma}(K)$ defined as follows. If $\sigma=\phi$, then $S_{\phi}(K)=K$. If $|\sigma|=0$, then $S_{\sigma}(K)$ is simplicial isomorphic to $K$ by replacing the unique vertex $w$ of $\sigma$ by a new vertex $v$. If $|\sigma|>0$, then the vertex set of $S_{\sigma}(K)$ is the vertex set of $K$ added with a new vertex $v$ and $S_{\sigma}(K)=(K{-}\{\sigma\}{*}{\rm link}_K\sigma)\cup C(\partial\sigma{*}{\rm link}_K\sigma)=\{\tau\!\in\!K\,|\,\sigma\!\not\subset\!\tau\}\cup
\{\{v\}{\cup}\sigma'{\cup}\tau\,|\,\sigma'\!\in\!\partial\sigma,\,\tau\!\in\!{\rm link}_K\sigma
\}$.
Two simplicial complexes are stellar equivalent if one can be obtained from another by a finite composite of stellar subdivision or its inverse. Two polyhedra are homeomorphic if and only if they have stellar equivalent triangulations(see \cite{a} and \cite{b}).
\vspace{3mm}

The following definition is a brief review of simplicial homology and cohomology theory. We always regard the simplicial chain complex $(C_*(K),d)$ and its dual cochain complex $(C^*(K),\delta)$ as the same chain group $C_*(K)=C^*(K)$ with different differentials $d$ and $\delta$. This trick is essential for the definition of (co)homeology spectral sequence.\vspace{1.5mm}

{\bf Definition 1.2} Let $K$ be a simplicial complex with vertex set $S$.  The $k$-th chain group $C_k(K)=C^k(K)$ of $K$ is the Abellian group generated by all ordered sequence $[v_0,v_1,\cdots,v_k]$ of elements of $S$ modular the following zero relations. $[v_0,\cdots,v_k]=0$ if $v_i=v_j$ for some $i\neq j$. If there is no repetition in $[v_0,\cdots,v_k]$, then $[v_0,\cdots,v_k]=0$ if $\{v_0,{\cdots},v_k\}\not\in K$ and $[\cdots,v_{i},\cdots,v_{j},\cdots]=-[\cdots,v_{j},\cdots,v_{i},\cdots]$ (the omitted part unchanged) for all $i<j$. Denote the unique generator of $C_{-1}(K)$ by $[\,]$ or $[\phi]$. $[v_0,{\cdots},v_k]\neq 0$ is called a chain simplex of $K$ of dimension $k$. The graded group
$C_*(K)=\oplus_{k\geqslant 0}C_k(K)=C^*(K)$ and $\w C_*(K)=\oplus_{k\geqslant -1}C_k(K)=\w C^*(K)$ are respectively the chain group and augmented chain group of $K$.

The simplicial chain complex $(C_*(K),d)$ and simplicial cochain complex $(C^*(K),\delta)$ of $K$ are defined as follows. For any chain simplex $[v_0,\cdots,v_k]$ ($\hat v_i$ means canceling the symbol from the term),\\
\hspace*{40mm}$d[v_0,\cdots,v_k]={\Sigma}_{i=0}^k(-1)^i[v_0,\cdots,\hat v_{i},\cdots,v_k],\,\,d[v_0]=0$,\\
\hspace*{50mm}$\delta[v_0,\cdots,v_k]={\Sigma}_{v\in S}[v,v_0,\cdots,v_n].$\\
$H_*(K)=H_*(C_*(K),d)$ and $H^*(K)=H^*(C^*(K),\delta)$ are respectively the homology and cohomology of $K$. For an Abellian group $G$, $H_*(K;G)=H_*(C_*(K){\otimes}G,d)$ and $H^*(K;G)=H^*(C^*(K){\otimes}G,\delta)$ are respectively the homology and cohomology of $K$ over $G$.

The reduced simplicial chain complex $(\w C_*(K),\w d)$ and reduced simplicial cochain complex $(\w C^*(K),\w\delta)$ of $K$ are defined as follows. For $\sigma\in\w C_*(K)$ with $|\sigma|>0$, $\w d\sigma=d\sigma$, $\w d[v]=[\,]$ and $\w d[\,]=0$. For $\sigma\in\w C^*(K)$ with $|\sigma|\geqslant 0$, $\w \delta\sigma=\delta\sigma$ and $\w \delta[\,]=\sum_{v\in S}[v]$. $\w H_*(K)=H_*(\w C_*(K),\w d)$ and $\w H^*(K)=H^*(\w C^*(K),\w\delta)$ are respectively the reduced homology and cohomology of $K$. For an Abellian group $G$, $\w H_*(K;G)=H_*(\w C_*(K){\otimes}G,d)$ and $\w H^*(K;G)=H^*(\w C^*(K){\otimes}G,\delta)$ are respectively the reduced homology and cohomology of $K$ over $G$.

If $L$ is a simplicial subcomplex of $K$, the relative chain group $C_*(K,L)=C_*(K)/C_*(L)=C^*(K,L)$. The relative simplicial chain complex $(C_*(K,L),d)$ is the quotient complex $(C_*(K),d)/(C_*(L),d)$. The relative simplicial cochain complex $(C^*(K,L),\delta)$ is the natural cochain subcomplex of $(C^*(K),\delta)$ generated by those chain simplexes of $C^*(K)$ that is not a chain simplex of $C^*(L)$. We have short exact sequences\\
\hspace*{45mm}$0\to C_*(L)\to C_*(K)\to C_*(K,L)\to 0$\\
\hspace*{45mm}$0\to C^*(K,L)\to C^*(K)\to C^*(L)\to 0$\\
that induce long exact sequences\\
\hspace*{25mm}$\cdots\to H_n(L;G)\to H_n(K;G)\to H_n(K,L;G)\to H_{n-1}(L;G)\to\cdots$\\
\hspace*{25mm}$\cdots\to H^{n-1}(L;G)\to H^n(K,L;G)\to H^n(K;G)\to H^{n}(L;G)\to\cdots$
\vspace{3mm}

{\bf Convention} Notice the difference between simplex and chain simplex. A chain simplex is a symbol and a simplex is a set. A chain simplex corresponds to a unique simplex consisting of all the vertices of the chain simplex. If there is no confusion, we use the same greek letter $\sigma,\tau,\cdots$ to denote both the chain simplex and the simplex it corresponds to. For two simplexes $\sigma$ and $\tau$, $\sigma{\cap}\tau$ and $\sigma{\cup}\tau$ are naturally defined. But for two chain simplexes $\sigma$ and $\tau$, $\sigma{\cap}\tau$ and $\sigma{\cup}\tau$ are not defined. For two chain simplexes $\sigma=[v_0,{\cdots},v_m]$ and $\tau=[w_0,{\cdots},w_n]$ such that the corresponding simplex $\sigma{\cap}\tau=\phi$, denote $\sigma{\cup}\tau=[v_0,{\cdots},v_m,w_0,{\cdots},w_n]$ and $\tau{\cup}\sigma=[w_0,{\cdots},w_n,v_0,{\cdots},v_m]$.
\vspace{3mm}

{\bf Definition 1.3} For a simplicial complex $K$, the double complex $(T^{*,*}(K),\Delta)$
\[\begin{array}{llllllll}&&&&&&T^{3,3}&\cdots\\&&&&&&_d\downarrow&\cdots\\
&&&&T^{2,2}&\stackrel{\delta}{\to}&T^{2,3}&\cdots\\&&&&_d\downarrow&&_d\downarrow&\cdots\\
&&T^{1,1}&\stackrel{\delta}{\to}&T^{1,2}&\stackrel{\delta}{\to}&T^{1,3}&\cdots\\
&&_d\downarrow&&_d\downarrow&&_d\downarrow&\cdots\\
T^{0,0}&\stackrel{\delta}{\to}&T^{0,1}&\stackrel{\delta}{\to}&T^{0,2}&\stackrel{\delta}{\to}&T^{0,3}&\cdots
\end{array}\]
is defined as follows. $ T^{*,*}(K)$ is the subgroup of $ C_*(K){\otimes} C^*(K)$ generated by all tensor product of chain simplexes $\sigma{\otimes}\tau$ such that $\sigma\!\subset\!\tau$ and $\Delta\colon T^{s,t}\to T^{s-1,t}{\oplus}T^{s,t+1}$ is defined by $\Delta(\sigma{\otimes}\tau)=(d\sigma){\otimes}\tau+(-1)^{|\sigma|}\sigma{\otimes}(\delta\tau)$. In another word, $(T^{*,*}(K),\Delta)$ is the double subcomplex of $(C_*(K),d){\otimes}(C^*(K),\delta)$ generated by all tensor product of chain simplexes $\sigma{\otimes}\tau$ such that $\sigma{\subset}\tau$. The cohomeology spectral sequence $\{ H_r^{s,t}(K),\Delta_r\}$ $r=0,1,\cdots$ of $K$ is induced by the horizontal filtration $F_{n}=\oplus_{s\leqslant n} T^{s,*}$ with $\Delta_r\colon H_r^{s,t}\to H_r^{s-r-1,t-r}$. Notice that to make the notation simple, the cohomeology spectral sequence starts from $r=0$. For $r>0$, $H_r^{*,*}(K)$ is called the $r$-th cohomeology group of $K$. $H_1^{*,*}(K)$ is simply denoted by $ H^{*,*}(K)$ and is called the cohomeology group of $K$. For an Abellian group $G$, define $ T^{*,*}(K;G)= T^{*,*}(K){\otimes}G$ and we have similar definition of  cohomeology spectral sequence or $r$-th cohomeology group $ H_r^{s,t}(K;G)$ of $K$ over $G$.

Similarly, the augmented double complex $(\w T^{*,*}(K),\w \Delta)$
\[\begin{array}{llllllll}&&&&&&\w T^{2,2}&\cdots\\&&&&&&_{\w d}\downarrow&\cdots\\
&&&&\w T^{1,1}&\stackrel{\w \delta}{\to}&\w T^{1,2}&\cdots\\&&&&_{\w d}\downarrow&&_{\w d}\downarrow&\cdots\\
&&\w T^{0,0}&\stackrel{\w \delta}{\to}&\w T^{0,1}&\stackrel{\w \delta}{\to}&\w T^{0,2}&\cdots\\
&&_{\w d}\downarrow&&_{\w d}\downarrow&&_{\w d}\downarrow&\cdots\\
\w T^{-1,-1}&\stackrel{\w \delta}{\to}&\w T^{-1,0}&\stackrel{\w \delta}{\to}&\w T^{-1,1}&\stackrel{\w \delta}{\to}&\w T^{-1,2}&\cdots
\end{array}\]
is the double subcomplex of $(\w C_*(K),\w d){\otimes}(\w C_*(K),\w \delta)$ generated by all product of chain simplexes $\sigma{\otimes}\tau$ such that $\sigma\!\subset\!\tau$. The reduced cohomeology spectral sequence $\{\w H_r^{s,t}(K),\w \Delta_r\}$ of $K$ is induced by the horizontal filtration $\w F_{n}=\oplus_{s\leqslant n}\w T^{s,*}$ with $\w \Delta_r\colon\w H_r^{s,t}\to\w H_r^{s-r-1,t-r}$. For $r>0$, $\w H_r^{s,t}(K)$ is called the $r$-th reduced cohomeology group of $K$.  $\w H_1^{*,*}(K)$ is simply denoted by $\w H^{*,*}(K)$ and is called the reduced cohomeology group of $K$. For an Abellian group $G$, define $\w T^{*,*}(K;G)=\w T^{*,*}(K){\otimes}G$ and we have similar definition of reduced cohomeology spectral sequence or reduced $r$-th cohomeology group $\w H_r^{s,t}(K;G)$ of $K$ over $G$.

Dually, define $(T_{*,*}(K),D)=((T^{*,*}(K))^*,\Delta^*)=({\rm Hom}_{\Bbb Z}(T^{*,*}(K),\Bbb Z),\Delta^*)$. Then we may regard $T_{*,*}(K)$ and $T^{*,*}(K)$ as the same group with differential $D$ defined by
\begin{eqnarray*}&&D([v_0,\cdots,v_m]{\otimes}[v_0,\cdots,v_m,v_{m+1},\cdots,v_{m+n}])\\
&=&\Sigma_{j=1}^n(-1)^{m+j}[v_0,\cdots,v_m]{\otimes}[v_0,\cdots,v_m,v_{m+1},\cdots,\hat v_{m+j},\cdots,v_{m+n}]\\
&&+\Sigma_{j=1}^{n}[v_{m+j},v_0,\cdots,v_m]{\otimes}[v_0,\cdots,v_m,v_{m+1},\cdots,v_{m+n}]\end{eqnarray*}
The homeology spectral sequence $\{ H^r_{s,t}(K),D^r\}$ $r=0,1,\cdots$ of $K$ is induced by the horizontal filtration $F_{n}=\oplus_{s\geqslant n} T_{s,*}$ with $D^r\colon H^r_{s,t}\to H^r_{s+r+1,t+r}$. For $r>0$, $H_r^{*,*}(K)$ is called the $r$-th homeology group of $K$. $H^1_{*,*}(K)$ is simply denoted by $ H_{*,*}(K)$ and is called the homeology group of $K$. For an Abellian group $G$, define $ T_{*,*}(K;G)= T_{*,*}(K){\otimes}G$ and we have similar definition of  homeology spectral sequence or $r$-th homeology group $ H^r_{s,t}(K;G)$ of $K$ over $G$. The reduced homeology spectral sequence and groups are similarly defined.\vspace{3mm}

{\bf Convention} All the spectral sequence in this paper starts from $r=0$ so that all the (co)homeology spectral sequence is an homeomorphism invariant from $r\geqslant 1$. All the definitions, theorems and proofs have their duals. We always give the dual definitions and theorems but the proof, we only give the natural one and omit the dual proof.\vspace{3mm}

{\bf Theorem 1.4}\, Let $K$ be a simplicial complex and $G$ be an Abellian group.  The reduced cohomeology and homeology spectral sequence of $K$ over $G$ satisfy that\\
\hspace*{15mm}$\w H_0^{s,s+t+1}(K;G)=\oplus_{|\sigma|=s}\, \w H^{t}({\rm link}_K\sigma;G)$,\,\,
$(\w H_0^{*,*}(K;G),\w \Delta_0)=(\oplus_{\sigma\in K}\,\w H^{*}({\rm link}_K\sigma;G),\w \Delta_0)$,\\
\hspace*{15mm}$\w H^0_{s,s+t+1}(K;G)=\oplus_{|\sigma|=s}\, \w H_{t}({\rm link}_K\sigma;G)$,\,\,
$(\w H^0_{*,*}(K;G),\w D^0)=(\oplus_{\sigma\in K}\,\w H_{*}({\rm link}_K\sigma;G),\w D^0)$,\\
with differential on the right side defined as follows. Suppose $\sigma=\{v_0,{\cdots},v_n\}\in K$ and $\sigma_i=\sigma{-}\{v_i\}$. Then  $\delta_{v_i}\colon(C_*({\rm link}_K\sigma),\delta)\to(C_*({\rm link}_K\sigma_i),\delta)$ and $d^{v_i}\colon(C_*({\rm link}_K\sigma_i),d)\to(C_*({\rm link}_K\sigma),d)$ are dual homomorphisms defined by $\delta_{v_i}(\tau)=\tau{\cup}[v_i]$ for all chain simplex $\tau\in{\rm link}_K\sigma$.
Denote $x{\in}\w H^*({\rm link}_K\tau;G)$ by $[x]_{\tau}$ in $\oplus_{\sigma\in K}\,\w H^{*}({\rm link}_K\sigma;G)$, then $\w \Delta_0([x]_{\tau})=\sum_{v_i}[\delta_{v_i}^*(x)]_{\tau_i}$ with $\tau_i{\cup}\{v_i\}=\tau$. Denote $y{\in}\w H_*({\rm link}_K\tau;G)$ by $[y]_{\tau}$ in $\oplus_{\sigma\in K}\,\w H_{*}({\rm link}_K\sigma;G)$, then $\w D^0([y]_{\tau})=\sum_{v_i}[d^{v_i}_*(y)]_{\tau_i}$ with $\tau{\cup}\{v_i\}=\tau_i$.

The cohomeology and homeology spectral sequence of $K$ over $G$ satisfy that \\
\hspace*{15mm}$H_0^{s,s+t+1}(K;G)=\oplus_{|\sigma|=s}\, \w H^{t}({\rm link}_K\sigma;G)$,\,\,
$(H_0^{*,*}(K;G),\Delta_0)=(\oplus_{\sigma\in K,\sigma\neq\phi}\,\w H^{*}({\rm link}_K\sigma;G),\Delta_0)$,\\
\hspace*{15mm}$H^0_{s,s+t+1}(K;G)=\oplus_{|\sigma|=s}\, \w H_{t}({\rm link}_K\sigma;G)$,\,\,
$(H^0_{*,*}(K;G),D^0)=(\oplus_{\sigma\in K,\sigma\neq\phi}\,\w H_{*}({\rm link}_K\sigma;G),D^0)$,\\
with differential on the right side defined as follows. $\Delta_0([x]_{\sigma})=\w\Delta_0([x]_{\sigma})$ if $|\sigma|\!>\!0$ and $\Delta_0([x]_{\sigma})=0$ if $|\sigma|\!=\!0$ and $(H^0_{*,*}(K),D^0)$ is a double subcomplex of $(\w H^0_{*,*}(K),\w D^0)$.

For any $n\!\geqslant\!0$, there are exact sequences
\[0\to \w H^{0,n}(K;G)\to H^{0,n}(K;G)\to\w H^n(K;G)\to \w H^{-1,n}(K;G)\to 0,\]
\[0\to \w H_{-1,n}(K;G)\to\w H_{n}(K;G)\to H_{0,n}(K;G)\to \w H_{0,n}(K;G)\to 0,\]
and $H^{s,n}(K;G)=\w H^{s,n}(K;G)$, $H_{s,n}(K;G)=\w H_{s,n}(K;G)$ if $s\!>\!0$.\vspace{3mm}

Proof\, We only prove the cohomeology case. The homeology case is totally dual.

$\w T^{u,v}(K)$ is generated by product of chain simplexes $\sigma\!\otimes\!\sigma'$ such that $\sigma\!\subset\!\sigma'$, $|\sigma|\!=\!u$, $|\sigma'|\!=\!v$. Equivalently, $\w T^{s,s+t+1}(K)$ is generated by all product of chain simplexes $\sigma\!\otimes\!(\tau{\cup}\sigma)$ such that $\sigma{\cap}\tau=\phi$, $|\sigma|=s$, $|\tau|=t$. So the correspondence $\sigma{\otimes}(\tau{\cup}\sigma)\to[\tau]_{\sigma}$ induces a graded group isomorphism $\varrho\colon\w T^{*,*}(K)\to \oplus_{\sigma\in K}\w C^*({\rm link}_K\sigma)$. Define $\w \Delta$ on the latter group by  $\w \Delta([\tau]_{\sigma})=
\sum_{v\in\sigma}[\delta_{v}(\tau)]_{\sigma-\{v\}}+(-1)^{|\sigma|}[\w \delta(\tau)]_{\sigma}$. Suppose $\sigma=[v_0,{\cdots},v_n]$ and $\tau=[w_0,{\cdots},w_m]$. Then
\begin{eqnarray*}&&\varrho\w \Delta([v_0,{\cdots},v_n]{\otimes}[w_0,{\cdots},w_m,v_0,{\cdots},v_n])\\
&=&\varrho(\Sigma_{i=0}^n(-1)^{i}[v_0,\cdots,\hat v_i,\cdots,v_n]{\otimes}[w_0,{\cdots},w_m,v_0,{\cdots},v_n])\\
&&+\varrho(\Sigma_{\{w,w_0,\cdots,w_m\}\in{\rm link}_K\sigma}(-1)^{n}[v_0,\cdots,v_n]{\otimes}[w,w_0,\cdots,w_m,v_0,\cdots,v_n])\\
&=&\Sigma_{i=0}^n[w_0,\cdots,w_n,v_i]_{\{v_0,\cdots,\hat v_i,\cdots,v_n\}}+(-1)^n[\w \delta([w_0,\cdots,w_m])]_{\{v_0,\cdots,v_n\}}\\
&=&\w \Delta([w_0,\cdots,w_m]_{\{v_0,\cdots,v_n\}})\\
&=&\w \Delta\varrho([v_0,{\cdots},v_n]{\otimes}[w_0,{\cdots},w_m,v_0,{\cdots},v_n]).
\end{eqnarray*}
So $\varrho$ is a well-defined double complex isomorphism that induces a spectral sequence isomorphism with $\w H_0^{*,*}$ isomorphism as in the theorem.

Define trivial double complex $U^{*,*}(K;G)$ by $U^{-1,t}(K;G)=\w H_0^{-1,t}(K;G)=\w H^t(K;G)$ and $U^{s,t}(K;G)=0$ if $s\!\geqslant\!0$. Then $U^{*,*}$ is a subcomplex of $(\w H_0^{*,*},\w \Delta_0)$ and the quotient complex is just $(H_0^{*,*},\Delta_0)$. So from the long exact sequence induced by the short exact sequence $0\to U^{*,*}\to \w H_0^{*,*}\to H_0^{*,*}\to 0$ we get the relation between the reduced cohomeology group and cohomeology group. \vspace{3mm}

From the proof of the above theorem we know that there are two ways to describe the structure of $(T^{*,*}(K;G),\Delta)$. One way is to denote the generators by $\sigma{\otimes}\sigma'$ with $\Delta(\sigma{\otimes}\sigma')=(d\sigma){\otimes}\sigma'+(-1)^{|\sigma|}\sigma{\otimes}(\delta\sigma')$. The other is to denote the generators by $[\tau]_{\sigma}$ ($\tau\in{\rm link}_K\sigma$) with $\Delta([\tau]_{\sigma})=\sum_{v\in\sigma}[\delta_{v}(\tau)]_{\sigma-\{v\}}+
(-1)^{|\sigma|}[\delta(\tau)]_{\sigma}$. We have to use both in later proofs. But we only denote the generators of $(H_0^{*,*}(K;G),\Delta_0)$ in one way, i.e., $[x]_{\sigma}$ with $x\in\w H^*({\rm link}_K\sigma;G)$. The reduced case and the homeology case are the same.\vspace{3mm}

{\bf Remark} We can similarly define $U_{s,t}(K)=C_s(K){\otimes}C^t(K^c)$ generated by all $\sigma{\otimes}\tau$ with $\sigma\in K$ and $\tau\not\in K$ such that $\sigma\subset\tau$ and $\Delta(\sigma{\otimes}\tau)=(d\sigma){\otimes}\tau+(-1)^{|\sigma|}\sigma{\otimes}(\delta\tau)$. Then the $R$-spectral sequence of $K$ $\{R_r^{s,t}(K),\Delta^r\}$ induced by the vertical filtration $R^{n}=\oplus_{s\geqslant n}U_{*,s}$ satisfies $R_0^{*,*}(K)=\oplus_{\sigma\not\in K}\w H_*(K|_{\sigma})$, where $K|_{\sigma}=\{\tau{\in}K\,|\,\tau{\subset}\sigma\}$. But this spectral sequence is neither a homeomorphism or homotopy invariant of $K$ nor a homeomorphism or homotopy invariant of the Alexander dual $K^*$ of $K$.\vspace{3mm}

{\bf Definition 1.5} Let $L$ be a simplicial subcomplex of $K$. $(T_{*,*}(L),D)$ and $(\w T_{*,*}(L),\w D)$ are respectively double subcomplexes of $(T_{*,*}(K),D)$ and $(\w T_{*,*}(K),\w D)$. Define relative double complex $(T_{*,*}(K,L),D)$ and $(\w T_{*,*}(K,L),\w D)$ respectively to be the corresponding quotient double complex and so we have a short exact sequence of double complexes\\
\hspace*{45mm}$0\to T_{*,*}(L)\to T_{*,*}(K)\to T_{*,*}(K,L)\to 0,$\\
\hspace*{45mm}$0\to\w T_{*,*}(L)\to\w T_{*,*}(K)\to\w T_{*,*}(K,L)\to 0.$\\
The relative homeology spectral sequence $H_{*,*}^r(K,L;G)$ and the relative reduced homeology spectral sequence $\w H_{*,*}^r(K,L;G)$ are the spectral sequence induced by horizontal filtration.

Dually, there are relative double complexes $(T^{*,*}(K,L),\Delta)$ and $(\w T^{*,*}(K,L),\w\Delta)$ and short exact sequence of double complexes\\
\hspace*{45mm}$0\to T^{*,*}(K,L)\to T^{*,*}(K)\to T^{*,*}(L)\to 0,$\\
\hspace*{45mm}$0\to\w T^{*,*}(K,L)\to\w T^{*,*}(K)\to\w T^{*,*}(L)\to 0.$\\
The relative cohomeology spectral sequence $H^{*,*}_r(K,L;G)$ and the relative reduced cohomeology spectral sequence $\w H^{*,*}_r(K,L;G)$ are the spectral sequence induced by horizontal filtration.
\vspace{3mm}

Notice that we have long exact sequence of $H_0$, $\w H_0$ and $H^0$, $\w H^0$ groups from the short exact sequence of double complex and relative double complex but can not generally get a long exact sequence of (reduced) homeology and (reduced) cohomeology groups.\vspace{3mm}

{\bf Definition 1.6} For a simplicial complex $K$, a block complex $\frak B_K=\{b_{n,\alpha}\}$ of $K$ is a set of simplicial subcomplexes $b_{n,\alpha}$ (called an $n$-block) of $K$ that satisfies the following conditions.

1) Every $b_{n,\alpha}$ is a triangulation of disk $D^n$. Specifically, the empty simplicial complex $b_{-1,\phi}$ is the only $(-1)$-block of $\frak B_K$ that is regarded as a triangulation of $D^{-1}=\phi$.

2) $|K|=\sqcup_{\,b_{n,\alpha}\in \frak B_K}(|b_{n,\alpha}|-|\partial b_{n,\alpha}|)$, where $|\cdot|$ denotes the underlying set of the geometrical realization space and $\sqcup$ is the disjoint union of sets.

Specifically, $\frak B^0_K=\{2^{\sigma}\}_{\sigma\in K}$ is a block complex of $K$ which is called the trivial block complex of $K$.

$b_{m,\alpha}$ is called a face of $b_{n,\beta}$ if $b_{m,\alpha}\subset b_{n,\beta}$. $b_{m,\alpha}$ is called a proper face of $b_{n,\beta}$ if $b_{m,\alpha}\subset\partial b_{n,\beta}$.\vspace{3mm}

From the definition of block complex we have the following conclusions that is very important for later proofs.

For $n\geqslant 0$, every boundary of an $n$-block is the union of some $(n{-}1)$-blocks, i.e.,  for every $b_{n,\alpha}\in \frak B_K$, there are $b_{n-1,\alpha_1},\cdots,b_{n-1,\alpha_k}\in \frak B_K$ such that $\partial b_{n,\alpha}=\cup_{i=1}^k b_{n-1,\alpha_i}$. For any two different blocks $b_{n,\alpha}$ and $b_{n,\beta}$ of the same dimension, $b_{n,\alpha}{\cap}b_{n,\beta}$ is the union of their mutual proper faces. For two different blocks $b_{m,\alpha}$ and $b_{n,\beta}$ such that $m<n$, either $b_{m,\alpha}$ is a proper face of $b_{n,\beta}$ or $b_{m,\alpha}{\cap}b_{n,\beta}=\{\phi\}$. For $n>0$ and every $\sigma\in \partial b_{n,\alpha}$, link$_{b_{n,\alpha}}\sigma$ is a triangulation of $D^{n-|\sigma|-1}$.\vspace{3mm}

{\bf Definition 1.7} Let $\frak B_K$ be a block complex of a simplicial complex $K$. $C_k(\frak B_K)=C^k(\frak B_K)$ is the free Abellian generated by all $k$-blocks of $\frak B_K$.  The graded group
$C_*(\frak B_K)=\oplus_{k\geqslant 0}C_k(\frak B_K)=C^*(\frak B_K)$ and $\w C_*(\frak B_K)=\oplus_{k\geqslant -1}C_k(\frak B_K)=\w C^*(\frak B_K)$ are respectively the block chain group and augmented block chain group of $\frak B_K$.

The block chain complex $(C_*(\frak B_K),d)$ and block cochain complex $(C^*(\frak B_K),\delta)$ of $\frak B_K$ are defined as follows. Take an orietation on every block of $\frak B_K$. That is, for $n>0$, all the chain $n$-simplexes of a given $n$-block $b_{n,\alpha}$ are devided into two classes, the class of positive chain simplexes and the class of negative chain simplexes. All the chain simplexes of a block of dimension $0$ or $-1$ are positive. For two blocks $b_{n,\alpha}$ and $b_{n-1,\beta}$, the connecting coefficient $[\alpha;\beta]$ is defined as follows. $[\alpha;\beta]=0$ if $b_{n-1,\beta}$ is not a face of $b_{n,\alpha}$; $[\alpha;\beta]=1$ if there is a positive chain simplex $[v_0,v_1\cdots,v_n]$ of $b_{n,\alpha}$ such that $[v_1,\cdots,v_{n}]$ is a positive chain simplex of $b_{n-1,\beta}$; $[\alpha;\beta]={-}1$ if there is a positive chain simplex $[v_0,v_1,\cdots,v_n]$ of $b_{n,\alpha}$ such that $[v_1,\cdots,v_{n}]$ is a negative chain simplex of $b_{n-1,\beta}$. It is obvious that the connecting coefficient is independent of the choice of the chain simplex. Then\\
\hspace*{30mm}$d(b_{n,\alpha})=\sum_{b_{n-1,\beta}\in \frak B_K}[\alpha;\beta]b_{n-1,\beta}$ for $n>0$, $d(b_{0,\alpha})=0,$\\
\hspace*{30mm}$\delta(b_{n,\alpha})=\sum_{b_{n+1,\beta}\in \frak B_K}[\alpha;\beta]b_{n+1,\beta}$ for $n\geqslant 0$.

The augmented block chain complex $(\w C_*(\frak B_K),\w d)$ and augmented block cochain complex $(\w C^*(\frak B_K),\w \delta)$ of $\frak B_K$ are defined by that $\w d(b_{n,\alpha})=d(b_{n,\alpha})$ if $n>0$, $\w d(b_{0,\alpha})=b_{-1,\phi}$, $\w d(b_{-1,\phi})=0$ and  $\w \delta(b_{n,\alpha})=\delta(b_{n,\alpha})$ if $n\geqslant 0$, $\w \delta(b_{-1,\phi})=\sum b_{0,\alpha}$ with $b_{0,\alpha}$ taken over all $0$-blocks.\vspace{3mm}

{\bf Definition 1.8} For a block complex $\frak B_K$ of a simplicial complex $K$, the double complex $(T^{*,*}(\frak B_K),\Delta)$
is the double subcomplex of $(C_*(\frak B_K),d){\otimes}(C^*(\frak B_K),\delta)$ generated by all tensor product of blocks $b_{m,\alpha}{\otimes}b_{n,\beta}$ such that $b_{m,\alpha}{\subset}b_{n,\beta}$. The cohomeology spectral sequence $\{ H_r^{s,t}(\frak B_K),\Delta_r\}$ $r=0,1,\cdots$ of $\frak B_K$ is induced by the horizontal filtration $F_{n}=\oplus_{s\leqslant n} T^{s,*}$ with $\Delta_r\colon H_r^{s,t}\to H_r^{s-r-1,t-r}$. $H_r^{*,*}(\frak B_K)$ is called the $r$-th cohomeology group of $\frak B_K$. $H_1^{*,*}(\frak B_K)$ is simply denoted by $ H^{*,*}(\frak B_K)$ and is called the cohomeology group of $\frak B_K$. For an Abellian group $G$, define $ T^{*,*}(\frak B_K;G)= T^{*,*}(\frak B_K){\otimes}G$ and we have similar definition of  cohomeology spectral sequence or $r$-th cohomeology group $ H_r^{s,t}(\frak B_K;G)$ of $\frak B_K$ over $G$.

The augmented double complex $(\w T^{*,*}(\frak B_K),\w\Delta)$
is the double subcomplex of $(\w C_*(\frak B_K),\w d){\otimes}(\w C^*(\frak B_K),\w\delta)$ generated by all tensor product of blocks $b_{m,\alpha}{\otimes}b_{n,\beta}$ such that $b_{m,\alpha}{\subset}b_{n,\beta}$. The reduced cohomeology spectral sequence $\{\w H_r^{s,t}(\frak B_K),\w\Delta_r\}$ of $\frak B_K$ is induced by the horizontal filtration $\w F_{n}=\oplus_{s\leqslant n}\w T^{s,*}$ with $\w\Delta_r\colon\w H_r^{s,t}\to\w H_r^{s-r-1,t-r}$. $\w H_r^{s,t}(\frak B_K)$ is called the $r$-th reduced cohomeology group of $\frak B_K$.  $\w H_1^{*,*}(\frak B_K)$ is simply denoted by $\w H^{*,*}(\frak B_K)$ and is called the reduced cohomeology group of $\frak B_K$. $\w H_r^{s,t}(\frak B_K;G)$ is similarly defined.

Dually, the double complex $(T_{*,*}(\frak B_K),D)$ is defined as follows. $T_{*,*}(\frak B_K)=T^{*,*}(\frak B_K)$ is the same bigraded Abellian group but the differential is defined by
\begin{eqnarray*}&&D(b_{m,\alpha}{\otimes}b_{n,\beta})\\
&=&\Sigma\, [\alpha;\sigma]b_{m+1,\sigma}{\otimes}b_{n,\beta}+\Sigma\, (-1)^m[\beta;\tau]b_{m,\alpha}{\otimes}b_{n-1,\tau},\end{eqnarray*}
where the sum is taken over all blocks such that $b_{m,\alpha}\subset b_{m+1,\sigma}\subset b_{n,\beta}$ and $b_{m,\alpha}\subset b_{n-1,\tau}\subset b_{n,\beta}$. The homeology spectral sequence $\{H^r_{s,t}(\frak B_K),D^r\}$ $r=0,1,\cdots$ of $\frak B_K$ is induced by the horizontal filtration $F_{n}=\oplus_{s\geqslant n} T^{s,*}$ with $\Delta^r\colon H^r_{s,t}\to H^r_{s+r+1,t+r}$. $H^r_{*,*}(\frak B_K)$ is called the $r$-th homeology group of $\frak B_K$. $H^1_{*,*}(\frak B_K)$ is simply denoted by $ H_{*,*}(\frak B_K)$ and is called the homeology group of $\frak B_K$. For an Abellian group $G$, define $T_{*,*}(\frak B_K;G)= T_{*,*}(\frak B_K){\otimes}G$ and we have similar definition of homeology spectral sequence or $r$-th homeology group $H^r_{s,t}(\frak B_K;G)$ of $\frak B_K$ over $G$.

The augmented double complex $(\w T_{*,*}(\frak B_K),\w D)$ is the dual double complex of $(\w T^{*,*}(\frak B_K),\w\Delta)$ with $\w D$ the generalization of $D$ in that $(-1)$-block is allowed in the formula. The reduced homeology spectral sequence or $r$-th homeology group $\w H^r_{*,*}(\frak B_K;G)$ of $\frak B_K$ over $G$ is similarly defined.
\vspace{3mm}

{\bf Theorem 1.9}\, Let $\frak B_K$ be a block complex of a simplicial complex $K$ and $G$ be an Abellian group.
For $b_{m,\alpha}\in \frak B_K$, define $C_n(\frak B_{m,\alpha})=C^n(\frak B_{m,\alpha})$ to be the free Abellian group generated by all blocks $b_{n,\beta}\in\frak B$ such that $b_{m,\alpha}\subset b_{n,\beta}$ and $C^*(\frak B_{m,\alpha})=\oplus_{n\geqslant m}C^n(\frak B_{m,\alpha})$. Let $(C^*(\frak B_{m,\alpha}),\delta)$ be the cochain subcomplex of $(C^*(\frak B_K),\delta)$ and denote its dual complex by $(C_*(\frak B_{m,\alpha}),d)$ which is not a chain subcomplex of $(C_*(\frak B_K),d)$. Define $\w H^*(\frak B_{m,\alpha};G)=H^*(C^*(\frak B_{m,\alpha}){\otimes}G,\delta)$ and $\w H_*(\frak B_{m,\alpha};G)=H_*(C_*(\frak B_{m,\alpha}){\otimes}G,d)$. For $b_{m,\alpha}\subset b_{s,\sigma}$, $(C^*(\frak B_{s,\sigma}),\delta)$ is a cochain subcomplex of $(C^*(\frak B_{m,\alpha}),\delta)$ and denote the inclusion map by $i_{\sigma,\alpha}$. The quotient map $j^{\alpha,\sigma}\colon (C_*(\frak B_{m,\alpha}),d)\to (C_*(\frak B_{s,\sigma}),d)$ is defined by $j^{\alpha,\sigma}(b_{n,\beta})=b_{n,\beta}$ if $b_{s,\sigma}\subset b_{n,\beta}$ and $j^{\alpha,\sigma}(b_{n,\beta})=0$ if $b_{s,\sigma}\not\subset b_{n,\beta}$. $i_{\sigma,\alpha}$ and $j^{\alpha,\sigma}$ are dual to each other.

The reduced cohomeology and homeology spectral sequence of $\frak B_K$ over $G$ satisfy that\\
\hspace*{7mm}$\w H_0^{m,n}(\frak B_K;G)=\oplus_{b_{m,\alpha}\in \frak B_K}\,\w H^n(\frak B_{m,\alpha};G)$,
$(\w H_0^{*,*}(\frak B_K;G),\w\Delta_0)=(\oplus_{b_{m,\alpha}\in \frak B_K}\,\w H^*(\frak B_{m,\alpha};G),\w\Delta_0)$,\\
\hspace*{7mm}$\w H^0_{m,n}(\frak B_K;G)=\oplus_{b_{m,\alpha}\in \frak B_K}\,\w H_n(\frak B_{m,\alpha};G)$,
$(\w H^0_{*,*}(\frak B_K;G),\w D^0)=(\oplus_{b_{m,\alpha}\in \frak B_K}\,\,\,\w H_*(\frak B_{m,\alpha};G),\w D^0)$,\\
where differentials are defined as follows. Denote $x\in\w H^*(\frak B_{m,\alpha};G)$ by $[x]_{\alpha}\in\oplus_{b_{m,\mu}\in \frak B_K}\,\w H^*(\frak B_{m,\mu};G)$, then $\w\Delta_0([x]_{\alpha})=\sum_{\alpha_i}[\alpha;\alpha_i][i^*_{\alpha,\alpha_i}(x)]_{\alpha_i}$, where $\alpha_i$ is taken over all faces $b_{m-1,\alpha_i}$ of $b_{m,\alpha}$. Denote $y\in\w H_*(\frak B_{m,\alpha};G)$ by $[y]_{\alpha}\in\oplus_{b_{m,\mu}\in \frak B_K}\,\w H_*(\frak B_{m,\mu};G)$, then $\w D^0([y]_{\alpha})=\sum_{\beta_j}[\alpha;\beta_j][j_*^{\alpha,\beta_j}(y)]_{\beta_j}$, where $\beta_j$ is taken over all blocks $b_{m+1,\beta_j}$ such that $b_{m,\alpha}\subset b_{m+1,\beta_j}$.

The cohomeology and homeology spectral sequence of $\frak B_K$ over $G$ satisfies that \\
\hspace*{7mm}$H_0^{m,n}(\frak B_K;G)=\oplus_{b_{m,\alpha}\in \frak B_K}\,\w H^n(\frak B_{m,\alpha};G)$,
$(H_0^{*,*}(\frak B_K;G),\Delta_0)=(\oplus_{b_{m,\alpha}\in \frak B_K,m>-1}\,\w H^*(\frak B_{m,\alpha};G),\Delta_0)$,\\
\hspace*{7mm}$H^0_{m,n}(\frak B_K;G)=\oplus_{b_{m,\alpha}\in \frak B_K}\,\w H_n(\frak B_{m,\alpha};G)$,
$(H^0_{*,*}(\frak B_K;G), D^0)=(\oplus_{b_{m,\alpha}\in \frak B_K,m>-1}\,\,\,\w H_*(\frak B_{m,\alpha};G),D^0)$,\\
where differentials are just the restriction differentials of the reduced case such that $b_{-1,\phi}$ is not allowed in the formula.

For any $n\!\geqslant\!0$, there are exact sequences\\
\hspace*{20mm}$0\to \w H^{0,n}(\frak B_K;G)\to H^{0,n}(\frak B_K;G)\to\w H^n(K;G)\to \w H^{-1,n}(\frak B_K;G)\to 0$\\
\hspace*{20mm}$0\to \w H_{-1,n}(\frak B_K;G)\to\w H_{n}(\frak B_K;G)\to\, H_{0,n}(K;G)\to\, \w H_{0,n}(\frak B_K;G)\to 0$\\
and $H^{s,n}(\frak B_K;G)=\w H^{s,n}(\frak B_K;G)$, $H_{s,n}(\frak B_K;G)=\w H_{s,n}(\frak B_K;G)$ if $s\!>\!0$.\vspace{2mm}

Proof\, By definition of horizontal filtration just as Theorem 1.4. But the proof here is more direct and does not need the isomorphism $\varrho$ in the proof of Theorem 1.4.
\vspace{3mm}

{\bf Theorem 1.10} Let $\frak B_K$ be a block complex of a simplicial complex $K$. Then for all Abellian group $G$ and all integers $r,s,t$ with $r>0$,\\
\hspace*{40mm}$H_r^{s,t}(\frak B_K;G)=H_r^{s,t}(K;G)$,\, $\w H_r^{s,t}(\frak B_K;G)=\w H_r^{s,t}(K;G)$,\\
\hspace*{40mm}$H^r_{s,t}(\frak B_K;G)=H^r_{s,t}(K;G)$,\, $\w H^r_{s,t}(\frak B_K;G)=\w H^r_{s,t}(K;G)$.\\
The cohomeology spectral sequence $H_r^{s,t}(K;G)$ and the homeology spectral sequence $H^r_{s,t}(K;G)$ respectively converge to the cohomology group $H^*(K;G)$ and the homology group $H_*(K;G)$. The reduced cohomeology spectral sequence $\w H_r^{s,t}(K;G)$ and the reduced homeology spectral sequence $\w H^r_{s,t}(K;G)$ converges to $G$ (at degree $0$). $H_r^{s,t}(K;G)$, $\w H_r^{s,t}(K;G)$, $H^r_{s,t}(K;G)$, $\w H^r_{s,t}(K;G)$ ($r>0$) are all homeomorphism invariants of the polyhedron $|K|$. Thus, for a polyhedron $X$, $H_r^{s,t}(X;G)$, $\w H_r^{s,t}(X;G)$, $H^r_{s,t}(X;G)$, $\w H^r_{s,t}(X;G)$ ($r>0$) are defined to be the corresponding group of any of its triangulation.\vspace{2mm}

Proof\, We only prove the unreduced cohomeology case . Other cases are similar.

Define $(T^{*,*}_{\frak B}(K),\Delta)$ to be the double subcomplex of $(C_*(K),d)\otimes(C^*(\frak B_K),\delta)$ generated by all $\sigma\otimes b_{m,\alpha}$ such that $\sigma\in b_{m,\alpha}$ and the group homomorphism $j\colon T^{*,*}(K)\to T^{*,*}_{\frak B}(K)$ by $j(\sigma{\otimes}\tau)=\sigma{\otimes}b_{m,\alpha}$ if $\tau\in b_{m,\alpha}{-}\partial b_{m,\alpha}$ and $|\tau|=m$, $j(\sigma{\otimes}\tau)=0$ if $\tau\in b_{m,\alpha}{-}\partial b_{m,\alpha}$ and $|\tau|<m$. It is easy to check that $j$ is a double complex homomorphism. The horizontal filtration also induces a cohomeology spectral sequence $H_r^{*,*}(T_{\frak B}^{*,*}(K);G)$ and so $j$ induces a spectral sequence homomorphism $j_r\colon H_r^{*,*}(K;G)\to H_r^{*,*}(T^{*,*}_{\frak B}(K);G)$. Consider\\
\hspace*{10mm}$j_0\colon H_0^{*,*}(K;G)=\oplus_{\sigma\in K-\{\phi\}}H^*({\rm link}_K\sigma;G)\to
H_0^{*,*}(T^{*,*}_{\frak B}(K);G)=\oplus_{\sigma\in K-\{\phi\}}H^*(C^*(B_{\sigma}){\otimes}G,\delta)$,\\
where $(C^*(B_{\sigma}),\delta)$ is the cochain subcomplex of $(C^*(\frak B_K),\delta)$ generated by
$\frak B_{\sigma}=\{b_{m,\alpha}\,|\,\sigma\in b_{m,\alpha}\}$. For $b_{m,\alpha}\in\frak B_{\sigma}$, define
$b^{\sigma}_{m,\alpha}=\{\tau\in{\rm link}_K\sigma\,|\,\sigma{\cup}\tau\in b_{m,\alpha}\}=
{\rm link}_{b_{m,\alpha}}\sigma$. By definition, $b^{\sigma}_{m,\alpha}$ is a triangulation of disk $D^{m-|\sigma|-1}$ with $b^{\sigma}_{m,\alpha}{-}\partial b^{\sigma}_{m,\alpha}=\{\tau\in{\rm link}_K\sigma\,|\,\sigma{\cup}\tau\in b_{m,\alpha}{-}\partial b_{m,\alpha}\}$. So
$\frak B'_{\sigma}=\{b^{\sigma}_{m,\alpha}\}$
is a block complex of ${\rm link}_K\sigma$ which is in 1-1 correspondence with $\frak B_{\sigma}$.
There is a monomorphism $i_{\frak B'_{\sigma}}\colon (C_*(\frak B'_{\sigma}),d)\to (C_*({\rm link}_K\sigma),d)$ defined by that for every block $b_{m,\alpha}\in\frak B'_{\sigma}$, $i_{\frak B'_{\sigma}}(b_{m,\alpha})=\sum\sigma$ with $\sigma$ taken over all positive chain $m$-simplexes of $b_{m,\alpha}$. It is known that $i_{\frak B'_{\sigma}}$ induces a homology isomorphism. Dually, there is an acyclic cochain subcomplex $(J^*(\frak B'_{\sigma}),\delta)$ of $(C^*({\rm link}_K\sigma),\delta)$ such that the quotient cochain complex is just $(C^*(\frak B'_{\sigma}),\delta)$ with the quotient map $j_{\frak B'_{\sigma}}\colon (C^*({\rm link}_K\sigma),\delta)\to (C^*(\frak B'_{\sigma}),\delta)$ defined by that $j_{\frak B'_{\sigma}}(\tau)=b_{m,\alpha}$ if $\tau\in b^{\sigma}_{m,\alpha}{-}\partial b^{\sigma}_{m,\alpha}$ and $|\tau|=m$ and that $j_{\frak B'_{\sigma}}(\tau)=0$ if $\tau\in b^{\sigma}_{m,\alpha}{-}\partial b^{\sigma}_{m,\alpha}$ and $|\tau|<m$. $j_{\frak B'_{\sigma}}$ induces cohomology isomorphism. Since
$(j\varrho^{-1})|_{C^*({\rm link}_K\sigma)}=j_{\frak B'_{\sigma}}$ ($\varrho$ as in the proof of Theorem 1.4), $j_0$ is an isomorphism. So $j_r\colon H_r^{*,*}(K;G)\to H_r^{*,*}(T^{*,*}_{\frak B}(K);G)$ is an isomorphism for all $r\geqslant 0$.

Define monomorphism $\psi\colon (T^{*,*}(\frak B_K),\Delta)\to (T^{*,*}_{\frak B}(K),\Delta)$ by that $\psi(b_{m,\alpha}{\otimes}b_{n,\beta})=\sum\sigma{\otimes}b_{n,\beta}$ with $\sigma$ taken over all positive chain $m$-simplexes of $b_{m,\alpha}$. Denote the quotient double complex by $(U_{\frak B}^{*,*}(K),\Delta)$. From the short exact sequence of double complexes $0\to T^{*,*}(\frak B_K)\stackrel{\psi}{\to} T^{*,*}_{\frak B}(K)\stackrel{\omega}{\to} U^{*,*}_{\frak B}(K)\to 0$ we have that the cohomeology spectral sequences induced by the horizontal filtration of the three double complexes satisfy that  $0\to H_0^{*,*}(\frak B_K);G)\stackrel{\psi_0}{\to} H_0^{*,*}(T^{*,*}_{\frak B}(K);G)\stackrel{\omega_0}{\to} H_0^{*,*}(U^{*,*}_{\frak B}(K);G)\to 0$ is still exact. If we have $H_1^{*,*}(U^{*,*}_{\frak B}(K);G)=0$, then from the long exact sequence of $H_1$-term induced by the $H_0$-term we have that $\psi_1$ is an isomorphism and inductively, $\psi_r\colon H_r^{*,*}(\frak B_K;G)\to H_r^{*,*}(T^{*,*}_{\frak B}(K);G)$ is an isomorphism for all $r>0$. Thus, $(j_r)^{-1}\psi_r\colon H_r^{*,*}(\frak B_K;G)\to H_r^{*,*}(K;G)$ is an isomorphism for all $r>0$.

Now we prove $H_1^{*,*}(U^{*,*}_{\frak B}(K);G)=0$. Define $\frak B_K^{(n)}=\cup_{b_{m,\alpha}\in \frak B_K,m\leqslant n}b_{m,\alpha}$ and the filtrations on $(H_0^{*,*},\Delta_0)$ as follows.
$F_n(\frak C)$ is the subgroup of $H_0^{*,*}(\frak C)$ generated by all $[x]_{\sigma}$ such that $\sigma\in \frak B_K^{(n)}$ ($\sigma\subset \frak B_K^{(n)}$ if $\sigma$ is a block), where $\frak C$ represents $(\frak B_K;G)$, $(T^{*,*}_{\frak B}(K);G)$ or $(U^{*,*}_{\frak B}(K);G)$. So we get a spectral sequence $F_r^{n,s,t}(\frak C)$ converging to $H_1^{*,*}(\frak C)$. By definition, \\
\hspace*{27mm}$F_0^{n,n,t}(\frak B_K;G)=\oplus_{b_{n,\alpha}\in\frak B_K}\w H^t_{n,\alpha}$,\, $\w H^t_{n,\alpha}=\w H^t(\frak B_{n,\alpha};G)$,\\
\hspace*{27mm}$F_0^{n,s,t}(T_{\frak B}^{*,*}(K);G)=\oplus_{|\sigma|=s,\,\sigma\in b_{n,\alpha}{-}\partial b_{n,\alpha}}\,\w H^{t}_{\sigma}$,\, $\w H^{t}_{\sigma}=\w H^t(\frak B_{n,\alpha};G)$.\\
The differential $\o\Delta_0$ of $F_0^{*,*,*}$ satisfies that $\o\Delta_0([x]_{b_{n,\alpha}})=0$ and $\o\Delta_0([y]_{\sigma})=\sum_i[\sigma;\sigma_i][y]_{\sigma_i}$, where $\sigma\in b_{n,\alpha}{-}\partial b_{n,\alpha}$ and $\sigma_i$ is taken over all $(|\sigma|{-}1)$-faces of $\sigma$ such that $\sigma_i\in b_{n,\alpha}{-}\partial b_{n,\alpha}$. Thus,\\
\hspace*{24mm}$(F_0^{*,*,*}(\frak B_K;G),\o\Delta_0)=(\oplus_{b_{n,\alpha}\in\frak B_K,n>-1}\w H^*(\frak B_{n,\alpha};G),\o\Delta_0)$ ($\o\Delta_0=0$),\\
\hspace*{24mm}$(F_0^{*,*,*}(T_{\frak B}^{*,*}(K);G),\o\Delta_0)=\oplus_{b_{n,\alpha}\in\frak B_K,n>-1}\,(C_*(b_{n,\alpha},\partial b_{n,\alpha}),d)\otimes\w H^*(\frak B_{n,\alpha};G)$.\\
Since $H_*(b_{n,\alpha},\partial b_{n,\alpha})=\Bbb Z$, the homomorphism $(\psi_0)_0\colon (F_0^{*,*,*}(\frak B_K;G),\o\Delta_0)\to(F_0^{*,*,*}(T_{\frak B}^{*,*}(K);G),\o\Delta_0)$ induces a cohomology isomorphism of $F_1$-term. So from the long exact sequence of $F_1$-term induced by the $F_0$-term we have that $F_1^{*,*,*}(U^{*,*}_{\frak B}(K);G)=0$. So $F_r^{*,*,*}(U^{*,*}_{\frak B}(K);G)=0$ for all $r>0$ and $H_1^{*,*}(U^{*,*}_{\frak B}(K);G)=0$.

From the vertical filtration $E_n=\oplus_{t\geqslant n}T^{*,t}(K)$ of $T^{*,*}(K)$, we get a spectral sequence $(E^{s,t}_r(K),\partial_r)$ with $\partial_r\colon E_r^{s,t}\to E_r^{s+r,t+r+1}$. By definition, $(E_0^{*,*}(K){\otimes}G,\partial)=(\oplus_{\sigma\in K-\{\phi\}}H_*(2^{\sigma};G),\partial)=(C^*(K){\otimes}G,\delta)$. So $E_1^{0,*}(K;G)=H^*(K;G)$ and $E_1^{s,*}(K;G)=0$ if $s>0$ and the spectral sequence collapse from $r>1$. Thus, $H^{*}(T^{*,*}(K)\otimes G,\Delta)=H^*(K;G)$ and the cohomeology spectral sequence converges to $H^*(K;G)$.

From the vertical filtration $\w E_n=\oplus_{t\geqslant n}\w T^{*,t}(K)$ of $\w T^{*,*}(K)$, we get a spectral sequence $(\w E^{s,t}_r(K),\w \partial_r)$ with $\w \partial_r\colon \w E_r^{s,t}\to \w E_r^{s+r,t+r+1}$. By definition, $(\w E_0^{*,*}(K){\otimes}G,\w \partial)=(\oplus_{\sigma\in K}\w H_*(2^{\sigma};G),\w \partial)=H_{-1}(\{\phi\};G)=G$. So $\w E_1^{0,0}(K;G)=G$ and $\w E_1^{s,t}(K;G)=0$ otherwise and the spectral sequence collapse from $r>1$. Thus, $H^{*}(\w T^{*,*}(K)\otimes G,\Delta)=G$ and the reduced cohomeology spectral sequence converges to $G$.

To prove that $H_r^{s,t}$ are homeomorphism invariants, we need only prove $H_r^{s,t}(S_{\sigma}(K);G)=H^r_{s,t}(K;G)$ for any stellar subdivision $S_{\sigma}$ on $K$. Take a block complex $\frak B_{S_{\sigma}(K)}$ on $S_{\sigma}(K)$ ($|\sigma|>0$) as follows. The blocks of $\frak B_{S_{\sigma}(K)}$ is in 1-1 correspondence with the simplexes of $K$. Denote by $b_{\tau}$ the block corresponding to $\tau\in K$. Then $b_{\tau}=2^{\tau}$ if $\sigma$ is not a face of $\tau$ and $b_{\tau}=S_{\sigma}(2^{\tau})$ if $\sigma\subset\tau$. It is obvious that the 1-1 correspondence $b_{\tau}\to 2^{\tau}$ from $\frak B_{S_{\sigma}(K)}$ to the trivial block complex $\frak B^0_K$ induces double complex isomorphism $(T^{*,*}(\frak B_{S_{\sigma}(K)}),\Delta)\cong (T^{*,*}(\frak B^0_K),\Delta)$. So the cohomeology spectral sequence is a homeomorphism invariant.
 \vspace{3mm}

For two simplicial complexes $K$ and $L$, their join is the simplicial complex $K*L=\{\sigma{\sqcup}\tau \,|\,\sigma{\in}K,\,\tau{\in}L\}$ ($\sqcup$ is the disjoint union) and their disjoint union is the simplicial complex $K\sqcup L=\{\phi\}\sqcup(K{-}\{\phi\})\sqcup(L{-}\{\phi\})$. Their Cartesian product $K\times L$ is the simplicial complex defined as follows. Suppose the vertex set of $K$ is $S=\{v_1,\cdots,v_m\}$ and the vertex set of $L$ is $T=\{w_1,\cdots,w_n\}$. The vertex set of $K\times L$ is $S\times T$ and a $k$-simplex of $K\times L$ is of the form $\{(v_{i_1},w_{j_1}),\cdots,(v_{i_k},w_{j_k})\}$ such that $i_s\leqslant i_{s+1}$, $j_s\leqslant j_{s+1}$ for $s=1,\cdots,k{-}1$ and $\{v_{i_1},\cdots,v_{i_k}\}\in K$, $\{w_{i_1},\cdots,w_{i_k}\}\in L$ (canceling repetition). If $\frak B_K$ and $\frak B_L$ are respectively block complex of $K$ and $L$, then their product block complex $\frak B_{K\times L}$ of $K\times L$ is the unique block complex satisfying that $\frak B_{K\times L}{-}\{\phi\}=(\frak B_K{-}\{\phi\})\times(\frak B_L{-}\{\phi\})$ and $|(b_{m,\alpha},b_{n,\beta})|=|b_{m,\alpha}|\times|b_{n,\beta}|$ for all $b_{m,\alpha}\in\frak B_K{-}\{\phi\}$ and $b_{n,\beta}\in\frak B_L{-}\{\phi\}$. \vspace{3mm}

{\bf Theorem 1.11} Let $K$, $L$ and $\frak B_K$, $\frak B_L$ be as above. Then
\[\begin{array}{ccc}
(T^{*,*}(K\sqcup L),\Delta)&=&(T^{*,*}(K),\Delta)\oplus (T^{*,*}(L),\Delta),\\ \\
\Sigma^{-1} (\w T^{*,*}(K{*}L),\Delta)&=&(\w T^{*,*}(K),\Delta)\,{\otimes}\,(\w T^{*,*}(L),\Delta),\\ \\
(T^{*,*}(\frak B_{K\times L}),\Delta)&=& (T^{*,*}(\frak B_K),\Delta)\,{\otimes}\,(T^{*,*}(\frak B_L),\Delta),\\\\
(T_{*,*}(K\sqcup L),D)&=&(T_{*,*}(K),D)\oplus (T_{*,*}(L),D),\\ \\
\Sigma^{-1} (\w T_{*,*}(K{*}L),D)&=&(\w T_{*,*}(K),D)\,{\otimes}\,(\w T_{*,*}(L),D),\\ \\
(T_{*,*}(\frak B_{K\times L}),D)&=& (T_{*,*}(\frak B_K),D)\,{\otimes}\,(T_{*,*}(\frak B_L),D),
\end{array}\]
where $\Sigma^{-1}$ means lowering the bidegree by $(-1,-1)$. The product isomorphism induces external product\\ \hspace*{40mm}$\times\colon H_r^{*,*}(K;G)\otimes H_r^{*,*}(L;H)\to H_r^{*,*}(K{\times}L;G{\otimes}H)$\\
\hspace*{40mm}$\times\colon H^r_{*,*}(K;G)\otimes H^r_{*,*}(L;H)\,\to H^r_{*,*}(K{\times}L;G{\otimes}H)$\\
for all Abellian groups $G,H$ and all $r\geqslant 0$. External products are bilinear.
\vspace{2mm}

Proof\, The union equality needs the conclusion ${\rm link}_{K*L}\,\sigma{*}\tau={\rm link}_K\sigma*{\rm link}_L\tau$. Everything else is by definition.\vspace{3mm}

{\bf Theorem 1.12} A simplicial complex $K$ is Cohen-Macaulay of dimension $n$ if for all simplexes $\sigma\in K$, $\w H^t({\rm link}_K\sigma)=0$ for all $t\neq n{-}|\sigma|{-}1$. For such a simplicial complex $K$, its reduced cohomeology spectral sequence over an Abellian group $G$ satisfies that for $r=1,\cdots,n{-}1$, the double complex $(\w H^{*,*}_r(K;G),\w \Delta_r)$ is as follows.
\[\begin{array}{ccccc}
   &&&&\w H^{n,n}_r=\w H^{n,n} \\
   & & &  &\vdots \\
   & & &  &\w  H^{r,n}_r=\w H^{r,n} \\
   & & & \stackrel{\w \Delta_r}{\swarrow} &  \\
   \w H^{-1,0}_r=\w H^{-1,0}_0 & \cdots & \w H^{-1,n-r}_r=\w H^{-1,n-r}_0 & &
 \end{array}\]
where all other $\w H^{s,t}_r(K;G)=0$ and $\w \Delta_r\colon \w H_{r}^{r,n}(K;G)\to\w H^{n-r}(K;G){=}\w H_r^{-1,n-r}(K;G)$ is an isomorphism. $\w \Delta_{n}\colon\w H_n^{n,n}(K;G)\to\w H_n^{-1,0}(K;G)$ is an epimorphism with kernel $G$ and $\w H_{n+1}^{n,n}(K;G)=G$, $\w H_{n+1}^{s,t}(K;G)=0$ otherwise. The reduced homeology spectral sequence satisfies  that for $r=1,\cdots,n{-}1$, the double complex $(\w H_{*,*}^r(K;G),\w D^r)$ is as follows.\\
\[\begin{array}{ccccc}
   &&&&\w H_{n,n}^r=\w H_{n,n} \\
   & & &  &\vdots \\
   & & &  &\w  H_{r,n}^r=\w H_{r,n} \\
   & & & \stackrel{\w D^r}{\nearrow} &  \\
   \w H_{-1,0}^r=\w H_{-1,0}^0 & \cdots & \w H_{-1,n-r}^r=\w H_{-1,n-r}^0 & &
 \end{array}\]
where all other $\w H_{s,t}^r(K;G)=0$ and $\w D^r\colon \w H_{n-r}(K;G){=}\w H^r_{-1,n-r}(K;G)\to\w H^{r}_{r,n}(K;G)$ is an isomorphism. $\w D^{n}\colon\w H^n_{-1,0}(K;G)\to\w H^n_{n,n}(K;G)$ is an monomorphism with cokernel $G$ and $\w H^{n+1}_{n,n}(K;G)=G$, $\w H^{n+1}_{s,t}(K;G)=0$ otherwise.\vspace{2mm}

Proof\, By definition $\w H^{s,t}_0(K;G)=0$ if $s\neq-1$ or $t\neq n$ and $\w \Delta_r\colon\w H_r^{s,t}\to \w H_r^{s-r-1,t-r}$. The conclusion is obvious.
\vspace{3mm}

{\bf Remark} The homomorphism $\w \Delta_r\colon\w H^{r,n}(K)\to\w H^{n-r}(K;G)$ and $\w D^r\colon\w H_{n-r}(K;G)\to\w H_{r,n}(K)$ for $r=1,\cdots,n$ in the above theorem is a generalization of Lefschetz duality theorem to Cohen-Macaulay complexes. For if $K$ is the triangulation of an orientable $n$-dimensional connected manifold $M$ with boundary, then $\w H_0^{-1,k}(M;G)=\w H^k(M;G)$ and we have the following commutative diagram of groups
\[\begin{array}{cccccccc}
\w H_0^{n,n}(K)&\stackrel{\w \Delta_0}{-\!\!-\!\!\!\to}&
\w H_0^{n-1,n}(K)&\stackrel{\w \Delta_0}{-\!\!-\!\!\!\to}\cdots \stackrel{\w \Delta_0}{-\!\!-\!\!\!\to}&
\w H_0^{0,n}(K)&\stackrel{\w \Delta_0}{-\!\!-\!\!\!\to}&
\w H_0^{-1,n}(K)=\Bbb Z&\to 0\\
 f_n\downarrow&&f_{n-1}\downarrow&&f_0\downarrow&&&\\
C_{n}(M,\partial M)&\stackrel{\w d}{-\!\!-\!\!\!\to}&
C_{n-1}(M,\partial M)&\stackrel{\w d}{-\!\!-\!\!\!\to}\cdots \stackrel{\w d}{-\!\!-\!\!\!\to}&
C_{0}(M,\partial M)&\to&
0&
\end{array}\]
where all $f_i$, $i>0$ is an isomorphism and $f_0$ is an epimorphism. So $\w H^{k,n}(K;G)=H_k(M,\partial M;G)$. Similarly, $\w H_{k,n}(K;G)=H^k(M,\partial M;G)$.\vspace{3mm}

In the following, we compute (co)homeology groups of some polyhedra.\vspace{3mm}

A $0$-dimensional simplicial complex $K_0$ is Cohen-Macaulay of dimension $0$, so $\w H_r^{0,0}(K_0;G)=G$ and $\w H^{s,t}_r(K_0;G)=0$ otherwise for all $r>0$. A pure $1$-dimensional simplicial complex $K_1$ (graphs with no isolated vertex) is Cohen-Macaulay of dimension $1$, so $\w H^{1,1}_r(K_1;G)=G$ and $\w H^{s,t}_r(K_1;G)=0$ otherwise for all $r>0$. Both the $n$-disk $D^n$ and the $n$-sphere $S^n$ are Cohen-Macaulay of dimension $n$, so $\w H^{n,n}_r(D^n;G)=\w H^{n,n}_r(S^n;G)=G$ and $\w H^{s,t}_r(D^n;G)=\w H^{s,t}(S^n;G)=0$ otherwise for all $r>0$. As an application, cohomeology group distinguishes disks of different dimension, i.e., $D^n\not\cong D^m$ if $m\neq n$.\vspace{3mm}

Let $C_nX$ be the join of $X$ with the $0$-dimensional simplicial complex with $n$ isolated vertices ($C_1X=CX$, $C_2X=SX$). By Theorem 1.11 and K\"{u}nneth Theorem, $\w H^{s+1,t+1}_r(C_nX;G)=\w H^{s,t}_r(X;G)$ for all $s,t\geqslant 0$ and $\w H^{s,t}_r(C_nX;G)=0$ if $s=-1$ or $t=-1$ for all $r>0$. Specifically, $\w H^{*,*}_r(C_3X;G)=\w H^{*,*}_r(SX;G)$. $C_3X$ and $SX\vee SX$ have the same homotopy type. But in general, $\w H^{*,*}(SX\vee SX;G)\neq\w H^{*,*}(C_3X;G)$. For example, if $X=S^n$, then $\w H^{n+1,n+1}(SX\vee SX)=\zz\oplus\zz$ and $\w H^{n+1,n+1}(C_3X)=\zz$. This shows that reduced cohomeology groups are not homotopy invariants.
\vspace{3mm}

A contractible polyhedron $X$ satisfies $H_r^{s,t}(X;G)=\w H_r^{s,t}(X;G)$ for all $r,s,t$ with $r>0$. Cohomeology group distinguishes contractible polyhedra by their dimension. Let $D^m\cup_{D^k}D^n$ be the contractible polyhedra defined as follows. Let $\phi_1\colon D^k\to D^m$ and $\phi_2\colon D^k\to D^n$ be the embedding of $D^k$ into the boundary of $D^m$ and $D^n$. $D^m\cup_{D^k}D^n$ is the quotient space $D^m{\sqcup} D^n/\!\!\sim$ with $\phi_1(x)\sim\phi_2(x)$ for all $x\in D^k$. If $k{+}1<{\rm min}(m,n)$, $H^{s,t}(D^m\cup_{D^k}D^n)=\zz$ for $(s,t)=(m,m),(n,n),(k,k{+}1)$ and $H^{s,t}(D^m\cup_{D^k}D^n)=0$ otherwise. If $k{+}1=m\leqslant n$, $H^{n,n}(D^m\cup_{D^k}D^n)=\zz$ and $H^{s,t}(D^m\cup_{D^k}D^n)=0$ otherwise.\vspace{3mm}

\section{Homeotopy and homeotopy type}\vspace{3mm}

{\bf Definition 2.1} A simplicial map $f\colon K\to L$ is solid if for all simplex $\sigma\in K$, $|f(\sigma)|=|\sigma|$.\vspace{3mm}

It is obvious that identity maps, inclusion maps from a simplicial subcomplex to the simplicial complex and composite maps of solid maps are all solid.\vspace{3mm}

{\bf Theorem 2.2} Let $f\colon K\to L$ be a solid map. Then $f$ induces double complex homomorphism\\
\hspace*{45mm}$f_*\colon T_{*,*}(K)\to T_{*,*}(L)$, \quad$\w f_*\colon \w T_{*,*}(K)\to \w T_{*,*}(L)$\\
and so homomorphisms\\
\hspace*{36mm}$f_*^r\colon H^r_{*,*}(K;G)\to H^r_{*,*}(L;G)$,\quad $\w f_*^r\colon \w H^r_{*,*}(K;G)\to \w H^r_{*,*}(L;G)$\\
for all Abellian group $G$ and $r\geqslant 0$ that satisfy the following.

1) The identity map $1_K$ induces identity isomorphisms $1\colon H^r_{*,*}(K;G)\to H^r_{*,*}(K;G)$ and $1\colon \w H^r_{*,*}(K;G)\to\w H^r_{*,*}(K;G)$.

2) The composite of solid maps induces the composite of induced homomorphisms, i.e., $(gf)_*^r=g_*^rf_*^r$ and $(\,\w g\,\w f\,)_*^r=\w g_*^r\,\w f_*^r$.

Dually, $f$ induces double complex homomorphism\\
\hspace*{45mm}$f^*\colon T^{*,*}(L)\to T^{*,*}(K)$,\quad $\w f^*\colon \w T^{*,*}(L)\to \w T^{*,*}(K)$\\
and so homomorphisms\\
\hspace*{36mm}$f^*_r\colon H_r^{*,*}(L;G)\to H_r^{*,*}(K;G)$,\quad$\w f^*_r\colon \w H_r^{*,*}(L;G)\to \w H_r^{*,*}(K;G)$\\
for all Abellian group $G$ and $r\geqslant 0$ that satisfy the following.

1) The identity map $1_K$ induces identity isomorphisms $1\colon H_r^{*,*}(K;G)\to H_r^{*,*}(K;G)$ and $1\colon \w H_r^{*,*}(K;G)\to\w H_r^{*,*}(K;G)$.

2) The composite of maps induces the adverse composite of induced homomorphisms, i.e., $(gf)^*_r=f^*_rg^*_r$ and $(\,\w g\,\w f\,)^*_r=\w f^*_r\,\w g^*_r$.

For a commutative ring $R$ and a fixed order on the vertex set of $K$, there is a cup product\\
\hspace*{50mm}$\cup\colon H^{*,*}_r(K;R){\otimes}H^{*,*}_r(K;R)\to H^{*,*}_r(K;R)$\\  
that makes $H^{*,*}_r(K;R)$ a graded commutative algebra over $R$ for all $r\geqslant 0$, where  graded commutative means for all $x\in H^{s,t}_r$ and $y\in H^{s',t'}_r$, $x{\cup}y=(-1)^{(s+t)(s'+t')}y{\cup}x$. For an order-preserving solid map $f\colon K\to L$, the induced homomorphism $f^*_r\colon H^{*,*}_r(L;R)\to H^{*,*}_r(K;R)$ is an algebra homomorphism for all $r\geqslant 0$. The algebra structure of $H^{*,*}_r(K;R)$ and algebra homomorphism $f^*_r$ are independent of the order of the vertex set when $r>0$.
\vspace{2mm}

Proof\, For $[v_0,\cdots,v_m]{\otimes}[v_0,\cdots,v_m,v_{m+1},\cdots,v_{m+n}]\in T_{*,*}(K)$, define\\
\hspace*{25mm}$f_*([v_0,\cdots,v_m]{\otimes}[v_0,\cdots,v_{m+n}])=
[f(v_0),\cdots,f(v_m)]{\otimes}[f(v_0),\cdots,f(v_{m+n})].$\\
It is obvious that for solid map $f$, $Df_*=f_*D$ and $\w D\w f_*=\w f_*\w D$. Dually, $\Delta f^*=f^*\Delta $ and $\w \Delta \w f^*=\w f^*\w \Delta $.

The cup product is the composite
$\cup\colon H^{*,*}_r(K;R)\otimes H^{*,*}_r(K;R)\stackrel{\times}{\to}H^{*,*}_r(K{\times}K;R{\otimes}R)
\stackrel{\mu_*}{\to}H^{*,*}_r(K{\times}K;R)
\stackrel{\Delta^*}{\to}H^{*,*}_r(K;R)$, where $\mu$ is the product map of $R$ and $\Delta$ is the diagonal map that is solid. In fact, we have a formula for the cup product. Suppose the vertex set of $K$ is $S=\{v_1,\cdots,v_n\}$. Then any chain simplex can be written in a unique way $[v_{i_0},\cdots,v_{i_k}]$ such that $i_0<\cdots<i_k$. Such a chain simplex is called an ordered simplex.
Then for $[v_{i_{m_0}},\cdots,v_{i_{m_s}}]{\otimes}[v_{i_0},\cdots,v_{i_k}],
[v_{j_{n_0}},\cdots,v_{j_{n_t}}]{\otimes}[v_{j_0},\cdots,v_{j_l}]\in T^{*,*}(K)$ with all chain simplexes ordered,
\begin{eqnarray*}&&\,\,\,([v_{i_{m_0}},\cdots,v_{i_{m_s}}]{\otimes}[v_{i_0},\cdots,v_{i_k}])\cup
([v_{j_{n_0}},\cdots,v_{j_{n_t}}]{\otimes}[v_{j_0},\cdots,v_{j_l}])\\
&=&\left\{\begin{array}{ll}
\hspace{9mm}(-1)^{kt}[v_{i_{m_0}},\cdots,\hat v_{i_{m_{s}}},v_{j_{n_0}},\cdots,v_{j_{n_t}}]{\otimes}
[v_{i_0},\cdots,v_{i_k},v_{j_1},\cdots,v_{j_l}] & {\rm if}\, i_k=j_0,\,m_s=k ,\\
\hspace{9mm}(-1)^{kt}[v_{i_{m_0}},\cdots,v_{i_{m_{s}}},\hat v_{j_{n_0}},\cdots,v_{j_{n_t}}]{\otimes}
[v_{i_0},\cdots,v_{i_k},v_{j_1},\cdots,v_{j_l}] & {\rm if}\, i_k=j_0,\,n_0=0, \\
(-1)^{st+kt+kl}[v_{j_{n_0}},\cdots,\hat v_{j_{n_{t}}},v_{i_{m_0}},\cdots,v_{i_{m_s}}]{\otimes}
{[}v_{j_0},\cdots,v_{j_t},v_{i_1},\cdots,v_{i_s}{]} & {\rm if}\, j_l=\,i_0,\,n_t=l,\\
(-1)^{st+kt+kl}[v_{j_{n_0}},\cdots,v_{j_{n_{t}}},\hat v_{i_{m_0}},\cdots,v_{i_{m_s}}]{\otimes}
{[}v_{j_0},\cdots,v_{j_t},v_{i_1},\cdots,v_{i_s}{]} & {\rm if}\, j_l=\,i_0,\,m_0=0,\\
\hspace{21mm}\,\, 0 & {\rm otherwise}.
             \end{array}\right.
\end{eqnarray*}
It is easy to check that the cup product makes $(T^{*,*}(K),\Delta)$ a differential graded algebra that is natural with respect to order-preserving solid map. Other conclusions are easy direct checkings. \vspace{3mm}

{\bf Definition 2.3} Let $X$ and $Y$ be two polyhedra. Two maps $f_0,f_1\colon X\to Y$ are homeotopic, denoted by $f_0\approx f_1$, if they satisfy the following.

1) $f_i$ is the geometrical realization map of solid map $f_i\colon K_i\to L_i$, $i=0,1$. So $K_i$ is a triangulation of $X$ and $L_i$ is a triangulation of $Y$.

2) There is a triangulation $\w K$ of $I{\times}X$ such that $K_i$ is a triangulation of $\{i\}{\times}X$ and a simplicial subcomplex of $\w K$, $i=0,1$. Denote the inclusion map from $K_i$ to $\w K$ by $\phi_i$.

3) There is a triangulation $\w L$ of $I{\times}Y$ such that $L_i$ is a triangulation of $\{i\}{\times}Y$ and a simplicial subcomplex of $\w L$, $i=0,1$. Denote the inclusion map from $L_i$ to $\w L$ by $\psi_i$.

4) There is a solid map $H\colon \w K\to \w L$ such that $H\phi_i=\psi_if_i$, $i=0,1$. $H$ is called a homeotopy (map) from $f_0$ to $f_1$.
\vspace{3mm}

{\bf Theorem 2.4} Homeotopic relation is an equivalence relation. Homeology group $H_{*,*}^r(-;G)$, reduced homeology group $\w H_{*,*}^r(-;G)$, cohomeology group $H^{*,*}_r(-;G)$ and reduced cohomeology group $\w H^{*,*}_r(-;G)$ are homeotopic invariants for all Abellian group $G$, i.e., homeotopic maps have the same induced (reduced) (co)homeology group homomorphism. Cohomeology algebra $H_r^{*,*}(-;R)$ over a commutative ring $R$ is also a homeotopic invariant. \vspace{2mm}

Proof\, To prove $f\approx f$, we have to prove that for any two $f_0,f_1$ with the same geometrical realization map $f$, $f_0\approx f_1$. If $f_0=f_1$, the homeotopy is obvious. If $f_0\neq f_1$, we need only prove the case $f_1$ is obtained from $f_0$ as follows. Suppose $f_0\colon K_0\to L_0$ and $\sigma\in L_0$. Since $f_0$ is solid, all the simplexes $\tau_1,\cdots,\tau_n$ of $K_1$ such that $f(\tau_i)=\sigma$ satisfy that $|\tau_i|=|\sigma|$. Then $f_1\colon S_{\tau_1}(\cdots(S_{\tau_n}(K_0)\cdots)\to S_{\sigma}(L_0)$ is the only simplicial map that has the same geometrical realization with $f_0$. Let $I_1$ be the triangulation of $I$ with two vertices $0,1$ and one edge. Then $\w\tau_i=\{1\}{\times}\tau_i$ and $\w \sigma=\{1\}{\times}\sigma$ are all simplexes of $I_1{\times} K_0$. $H\colon S_{\w\tau_1}(\cdots(S_{\w\tau_n}(I_1{\times}K_0)\cdots)\to S_{\w\sigma}(I_1{\times}L_0)$ is defined by that $H(i,v)=f_0(v)$ for all vertex $v\in K_0$ and $i=0,1$ and $H(1,w_i)=v$ for all the new vertex $w_i$ of $S_{\tau_i}(K_0)$ and $v$ is the new vertex of $S_{\sigma}(L_0)$. $H$ is obvously a homeotopy from $f_0$ to $f_1$. Suppose $f_0\approx f_1$ with $K_i, L_i,\w K,\w L,H$ as in the definition. Define $K^{-1}_i=K_{1-i}$, $L^{-1}_i=L_{1-i}$, $\w K^{-1}$ is isomorphic to $\w K$ by the 1-1 correspondence $(t,x)\to (1{-}t,x)$ of $I{\times} X$,  $\w L^{-1}$ is isomorphic to $\w L$ by the 1-1 correspondence $(t,y)\to (1{-}t,y)$ of $I{\times} Y$, $H^{-1}(t,x)=H(1{-}t,x)$ (geometrical realization). Then $f_1\approx f_0$ with  $K^{-1}_i, L^{-1}_i,\w K^{-1},\w L^{-1},H^{-1}$ as in the definition.  Suppose $f_0\approx f_1$ with $K_i, L_i,\w K,\w L,H$ as in the definition, $i=0,1$ and $f_1\approx f_2$ with $K_i, L_i,\w K',\w L',H'$ as in the definition, $i=1,2$ and the vertex set of $\w K,\w K',\w L,\w L'$ are respectively $\w S,\w S',\w T,\w T'$. Define the triangulation $\w M$ of $I{\times}X$ and the triangulation $\w N$ of $I{\times}Y$ with respectively vertex set $U$ and $V$ by\\
\hspace*{30mm}$U=\{(\frac 12t,x)\,|\,(t,x)\in \w S\,\}\cup\{(\frac12+\frac 12t,x)\,|\,(t,x)\in \w S'\,\}$,\\
\hspace*{30mm}$V=\{(\frac 12t,x)\,|\,(t,x)\in \w T\,\}\cup\{(\frac12+\frac 12t,x)\,|\,(t,x)\in \w T'\}$.\\
Define $(H{*}H')(t,x)=H(2t,x)$ if $0\leqslant t\leqslant\frac 12$ and $(H{*}H')(t,x)=H'(2t-1,x)$ if $\frac 12\leqslant t\leqslant 1$. Then $f_0\approx f_2$ with $K_0,K_2,L_0,L_2,\w M,\w N, H{*}H'$ as in the definition. So homeotopic relation is an equivalence relation.

The conclusion that (co)homeology groups are homeotopic invariant is a corollary of the following theorem. The reduced (co)homeology groups case naturally follows from the unreduced case.
\vspace{3mm}

{\bf Theorem 2.5} Let $K,L,\w K,\w L$ be respectively triangulation of $X,Y,I{\times}X,I{\times}Y$. $\Phi\colon\w K\to\w L$ and $f\colon K\to L$ are solid maps such that $\Phi\phi=\psi f$, where $\phi\colon X{=}\{0\}{\times}X\to I{\times}X$ and $\psi\colon Y{=}\{0\}{\times}Y\to I{\times}Y$ are inclusion maps. Then the induced (co)homeology group (over any Abellian group) homomorphism of $\Phi$ and $f$ satisfy that $\Phi^*_r=\iota^*{\times} f^*_r$ and $\Phi_*^r=\iota\,{\times} f_*^r$ from $r>0$, where $\iota$ is the generator class of $H_{1,1}(I)=H^r_{1,1}(I)$ and $\iota^*$ is the dual class of $\iota$ and $\times$ is the external product $H^r_{*,*}(I)\otimes H^r_{*,*}(Y;G)\to H^r_{*,*}(I{\times}Y;G)$ and $H_r^{*,*}(I)\otimes H_r^{*,*}(X;G)\to H_r^{*,*}(I{\times}X;G)$.\vspace{2mm}

Proof\, Define $1$-dimensional simplicial complex $I_n$ ($n{>}0$) as follows. The vertex set is $\{v_i\,|\,i=0,1,{\cdots},n\}$, the edge set is $\{\{v_i,v_{i+1}\}\,|\,i=0,1,\cdots,n{-}1\}$. Then $\{I_n\}$ is the set of all triangulations of $I$. $H_{s,t}(I_n)=0$ if $(s,t)\neq(1,1)$ and $H_{1,1}(I_n)=\Bbb Z$ with generator class represented by $[\phi]_{\{v_i,v_{i+1}\}}$ for any $0\leqslant i<n$. This implies that any solid map from $I_m$ to $I_n$ induces identity isomorphism of homeology groups. So by Theorem 1.11 and K\"{u}nneth theorem, for any polyhedron $Z$ and solid map $g\colon I\to I$, the map $g{\times}1_Z\colon I{\times}Z\to I{\times}Z$ induces identity isomorphism of homeology groups $H^r_{*,*}(I{\times}Z;G)$ for all Abellian group $G$.

Define $\o K$ to be a triangulation of $I{\times}X$ such that $\w K$ is a simplicial subcomplex of $\o K$ with $|\w K|=[\frac 12,1]{\times}X$ and that $I_1{\times}K$ is a simplicial subcomplex of $\o K$ with $|I_1{\times}K|=[0,\frac 12]{\times}X$. $\o L$ is similarly defined. Then we have the following commutative diagram of simplicial maps.
\[\begin{array}{ccc}
\quad\w K & \stackrel{\Phi}{\longrightarrow} &\quad \w L \\
^{\tilde i}\downarrow &  & ^{\tilde j}\downarrow \\
\quad\o K & \stackrel{\o\Phi}{\longrightarrow} & \quad\o L \\
_i\uparrow &  & _j\uparrow \\
I_1{\times}K & \stackrel{1{\times} f}{\longrightarrow} & I_1{\times}L
  \end{array}\]
where $i,\w i,j,\w j$ are the inclusion maps. The geometrical realization map of the above diagram is
\[\begin{array}{ccc}
\quad\quad I{\times}X & \stackrel{\Phi}{\longrightarrow} &\quad\quad I{\times}Y \\
^{g_0\times 1_X}\downarrow &  &^{g_0\times 1_Y}\downarrow \\
\quad\quad I{\times}X & \stackrel{\o\Phi}{\longrightarrow} &\quad\quad I{\times}Y \\
_{g_1\times 1_X}\uparrow &  & _{g_1\times 1_Y}\uparrow \\
\quad\quad I{\times}X & \stackrel{1{\times} f}{\longrightarrow} &\quad\quad I{\times}Y
  \end{array}\]
where $g_0(t)=\frac 12(1{+}t)$ and $g_1(t)=\frac 12t$ for all $t\in I$. By the above conclusion, all the vertical maps induce identity isomorphisms. If simplicial maps that have the same geometrical realization map have the same induced homeology group homomorphism, then the homeology group homomorphisms induced by $i,\w i,j,\w j$ are all identity isomorphisms and so $\Phi_*^r=(1{\times}f)_*^r=\iota{\times}f_*^r$. Now we prove that simplicial maps that have the same geometrical realization map have the same induced homeology group homomorphism. We need only prove the case $f_1$ is obtained from $f_0$ as in the proof of $f\approx f$ of Theorem 2.4. This case is obvious.
\vspace{3mm}

{\bf Definition 2.6} Two polyhedra $X$ and $Y$ are homeotopic equivalent, or have the same homeotopy type, denoted by $X\approxeq Y$, if there are geometrical realization maps of solid maps $f\colon X\to Y$ and $g\colon Y\to X$ such that $gf\approx 1_X$ and $fg\approx 1_Y$. $f$ is called a homeotopy equivalence from $X$ to $Y$.\vspace{3mm}

Homeotopy type is a coarser relation than homeomorphism. For example, the following two graphs are of the same homeotopy type but not homeomorphic.\vspace{3mm}

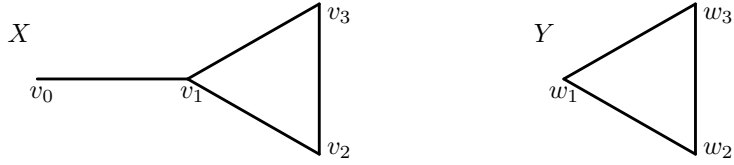
\begin{center}
\setlength{\unitlength}{0.5mm}
\begin{picture}(200,40)(0,0)
\Thicklines
\drawline(0,20)(40,20)
\drawline(40,20)(75,0)
\drawline(40,20)(75,40)
\drawline(75,0)(75,40)

\drawline(140,20)(175,0)
\drawline(140,20)(175,40)
\drawline(175,0)(175,40)

\put(-8,30){{$X$}}
\put(132,30){{$Y$}}
\put(-2,15){{$v_0$}}
\put(38,15){{$v_1$}}
\put(77,0){{$v_2$}}
\put(77,36){{$v_3$}}
\put(136,15){{$w_1$}}
\put(177,0){{$w_2$}}
\put(177,36){{$w_3$}}
\end{picture}
\end{center}\vspace{3mm}

$f\colon X\to Y$ is defined by $f(v_i)=w_i$ for $i=1,2,3$ and $f(v_0)=w_3$. $g\colon Y\to X$ is the inclusion $g(w_i)=v_i$ for $i=1,2,3$. The homeotopy from $1_X$ to $gf$ is shown in the following picture, where the homeotopy extends the identity of the triangle cylinder by maping the left two simplexes $1$ and $2$ respectively to the right two simplexes $1$ and $2$.\vspace{3mm}

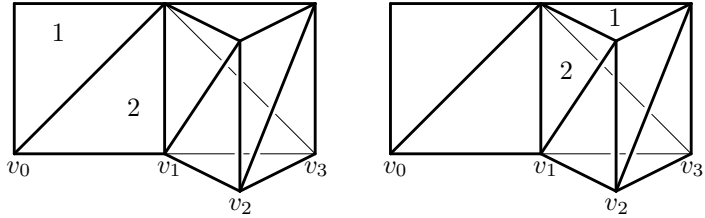
\begin{center}
\setlength{\unitlength}{0.5mm}
\begin{picture}(200,40)(0,0)
\Thicklines
\drawline(0,0)(0,40)
\drawline(0,0)(40,0)
\drawline(0,40)(40,40)
\drawline(40,0)(40,40)
\drawline(0,0)(40,40)
\drawline(40,0)(60,-10)
\drawline(60,-10)(80,0)
\drawline(40,40)(60,30)
\drawline(60,30)(80,40)
\drawline(80,0)(80,40)
\drawline(40,0)(60,30)
\drawline(60,-10)(80,40)
\drawline(40,40)(80,40)
\drawline(60,-10)(60,30)
\thinlines
\drawline(40,0)(58.5,0)
\drawline(61,0)(62.5,0)
\drawline(65.5,0)(80,0)
\drawline(40,40)(55,25)
\drawline(57,23)(59,21)
\drawline(61,19)(67.5,12.5)
\drawline(69.5,10.5)(80,0)

\put(10,30){1}
\put(30,10){2}

\Thicklines
\drawline(100,0)(100,40)
\drawline(100,0)(140,0)
\drawline(100,40)(140,40)
\drawline(140,0)(140,40)
\drawline(100,0)(140,40)
\drawline(140,0)(160,-10)
\drawline(160,-10)(180,0)
\drawline(140,40)(160,30)
\drawline(160,30)(180,40)
\drawline(180,0)(180,40)
\drawline(140,0)(160,30)
\drawline(160,-10)(180,40)
\drawline(140,40)(180,40)
\drawline(160,-10)(160,30)
\thinlines
\drawline(140,0)(158.5,0)
\drawline(161,0)(162.5,0)
\drawline(165.5,0)(180,0)
\drawline(140,40)(155,25)
\drawline(157,23)(159,21)
\drawline(161,19)(167.5,12.5)
\drawline(169.5,10.5)(180,0)

\put(158,34){1}
\put(145,20){2}

\put(-2,-5){{$v_0$}}
\put(38,-5){{$v_1$}}
\put(57,-15){{$v_2$}}
\put(77,-5){{$v_3$}}

\put(98,-5){{$v_0$}}
\put(138,-5){{$v_1$}}
\put(157,-15){{$v_2$}}
\put(177,-5){{$v_3$}}
\end{picture}
\end{center}
\vspace{7.5mm}

Homeotopy type is finer than homotopy type. For example, $D^n$ and $D^m$ ($m\neq n$) are not of the same homeotopy type since they have different (co)homeology groups. But they are of the same homotopy type. Another example is $C_3S^n$ and $S^{n+1}\vee S^{n+1}$ (see the examples in the end of the previous section). Just like that (co)homology group is the convergence group of (co)homeology spectral sequence, homotopy relation is the limit of homeotopy relation in the sense of the following theorem.\vspace{2mm}

{\bf Theorem 2.7} For a map $f\colon X\to Y$, its graph is the map $\w f\colon X\to X{\times}Y$ defined by $\w f(x)=(x,f(x))$ for all $x\in X$. If $f$ is the gemeotrical realization of a simplicial map, then $\w f$ is a solid map. If $f,g$ are both gemeotrical realization maps, then $f\sim g$ if and only if $\w f\approx\w g$.\vspace{2mm}

Proof\, Suppose $f$ is the geometrical realization of $f\colon K\to L$. The vertex set of $L$ is $\{w_1,\cdots,w_n,\cdots\}$ and the vertex set of $K$ is $\{v_1,\cdots,v_m\}$ such that there are $1\leqslant i_1<i_2<\cdots<i_{n-1}<i_n=m$ with $f(v_{j})=w_k$ for $j=i_{k-1}{+}1,\cdots,i_k$ ($i_0=0$). $f$ is an order preserving map and so $\w f\colon K\to K{\times}L$ is a solid map that maps the ordered simplex $\{v_{u_1},\cdots,v_{u_s}\}$ of $K$ to the ordered simplex $\{(v_{u_1},f(v_{u_1})),\cdots,(v_{u_s},f(v_{u_s}))\}$ of $K{\times}L$.

Suppose $f\sim g$, then there is geometrical realization homotopy $h\colon I{\times}X\to Y$ from $f$ to $g$. The graph $\w h\colon I{\times}X\to I{\times}X{\times}Y$ is just the homeotopy from $\w f$ to $\w g$. Conversely, if $\w f\approx \w g$ and $\w h\colon I{\times}X\to I{\times}X{\times}Y$ is the homeotopy from $\w f$ to $\w g$, then $h=p\,\w h$ is the homotopy from $f$ to $g$, where $p\colon I{\times}X{\times}Y\to Y$ is the projection map.
\vspace{3mm}

{\bf Theorem 2.8} Let $f\colon X\to Y$ be a geometrical realization map and $G$ be a module over the commutative ring $R$, then $f$ induces on the bigraded group $H^{*,*}_r(X;G)$ a bigraded module structure over the  bigraded algebra $H^{*,*}_r(Y;R)$ for all $r>0$. If $f\sim g$, then they induce the same module structure on $H^{*,*}_r(X;G)$.\vspace{2mm}

Proof\, The product map of the module is the composite of the following
\[H^{*,*}_r(X;G)\otimes H^{*,*}_r(Y;R)\stackrel{\times}{\longrightarrow}
H^{*,*}_r(X{\times}Y;G{\otimes}R)\stackrel{\mu_*}{\longrightarrow}H^{*,*}_r(X{\times}Y;G)
\stackrel{(\w f)^*}{\longrightarrow}H^{*,*}_r(X;G),\]
where $\mu$ is the module map of $G$. By Theorem 2.7, homotopic maps induce the same module structure.
\vspace{3mm}

This theorem shows that for any continuous map (need not be a geometrical realization map!) $f$ from polyhedron $X$ to polyhedron $Y$, the module structure of $H^{*,*}_r(X;G)$ over $H^{*,*}_r(Y;R)$ induced by $f$ can be defined.
By the simplicial approximation theorem, we need only define the module structure to be the one induced by any geometrical realization map homotopic to $f$.
\end{document}